\documentclass[12pt]{amsart}
\usepackage{amsmath, amssymb}
\usepackage{graphicx,color,hyperref}
\usepackage{verbatim, enumitem}

\newcommand\cL{{\mathcal L}}

\newcommand\bN{{\mathbb N}}
\newcommand\bC{{\mathbb C}}

\def\ii{{\rm i}\kern1pt}

\def\PSH{\mathop{\rm PSH}\nolimits}

\def\proof{\noindent{\it Proof.} }

\def\build#1^#2_#3{\mathrel{\mathop{\null#1}\limits^{#2}_{#3}}}
\def\mertorelbar{\vrule width0.6ex height0.65ex depth-0.55ex}
\def\merto{\mathrel{\mertorelbar\kern1.3pt\mertorelbar\kern1.3pt\mertorelbar
    \kern1.3pt\mertorelbar\kern-1ex\raise0.28ex\hbox{${\scriptscriptstyle>}$}}}

\catcode`\@=11
\newdimen\@rrowlength \@rrowlength=6ex
\def\ssrelbar{\vrule width\@rrowlength height0.64ex depth-0.56ex\kern-4pt}
\def\llra#1{\@rrowlength=#1\ssrelbar\rightarrow}
\catcode`\@=12

\def\semidirect{\mathop{\kern2pt\vrule depth-0.3pt height4.3pt 
\kern-2pt\times}\nolimits}

\newdimen\plainitemindent \plainitemindent=18pt
\def\plainitem#1{\par\noindent
\hangindent\plainitemindent\hbox to\plainitemindent{#1\hss}\ignorespaces}

\catcode`\@=11
\def\openup{\afterassignment\@penup\dimen@=}
\def\@penup{\advance\lineskip\dimen@
  \advance\baselineskip\dimen@
  \advance\lineskiplimit\dimen@}
\newdimen\jot \jot=3pt
\newskip\plaincentering \plaincentering=0pt plus 1000pt minus 1000pt
\def\ialign{\everycr{}\tabskip\z@skip\halign}
\def\eqalign#1{\null\,\vcenter{\openup\jot\m@th
  \ialign{\strut\hfil$\displaystyle{##}$&$\displaystyle{{}##}$\hfil
      \crcr#1\crcr}}\,}
\newif\ifdt@p
\def\displ@y{\global\dt@ptrue\openup\jot\m@th
  \everycr{\noalign{\ifdt@p \global\dt@pfalse \ifdim\prevdepth>-1000\p@
      \vskip-\lineskiplimit \vskip\normallineskiplimit \fi
      \else \penalty\interdisplaylinepenalty \fi}}}
\def\@lign{\tabskip\z@skip\everycr{}} 
\def\displaylines#1{\displ@y \tabskip\z@skip
  \halign{\hbox to\displaywidth{$\@lign\hfil\displaystyle##\hfil$}\crcr
    #1\crcr}}
\def\eqalignno#1{\displ@y \tabskip\plaincentering
  \halign to\displaywidth{\hfil$\@lign\displaystyle{##}$\tabskip\z@skip
    &$\@lign\displaystyle{{}##}$\hfil\tabskip\plaincentering
    &\llap{$\@lign##$}\tabskip\z@skip\crcr
    #1\crcr}}
\def\leqalignno#1{\displ@y \tabskip\plaincentering
  \halign to\displaywidth{\hfil$\@lign\displaystyle{##}$\tabskip\z@skip
    &$\@lign\displaystyle{{}##}$\hfil\tabskip\plaincentering
    &\kern-\displaywidth\rlap{$\@lign##$}\tabskip\displaywidth\crcr
    #1\crcr}}
\def\plaincases#1{\left\{\,\vcenter{\normalbaselines\m@th
    \ialign{$##\hfil$&\quad##\hfil\crcr#1\crcr}}\right.}
\def\plainmatrix#1{\null\,\vcenter{\normalbaselines\m@th
    \ialign{\hfil$##$\hfil&&\quad\hfil$##$\hfil\crcr
      \mathstrut\crcr\noalign{\kern-\baselineskip}
      #1\crcr\mathstrut\crcr\noalign{\kern-\baselineskip}}}\,}

\catcode`\@=12

\def\dlraw{\mathrel{\rlap{$\longrightarrow$}\kern-1pt\longrightarrow}}
\def\vlra{\mathrel{\smash-}\joinrel\mathrel{\smash-}\joinrel%
\kern-2pt\longrightarrow}
\def\srelbar{\vrule width0.6ex height0.65ex depth-0.55ex}
\def\merto{\mathrel{\srelbar\kern1.3pt\srelbar\kern1.3pt\srelbar
    \kern1.3pt\srelbar\kern-1ex\raise0.28ex\hbox{${\scriptscriptstyle>}$}}}

\newdimen\claimskip \claimskip=7pt
\long\def\claim#1|#2\endclaim
{\removelastskip\vskip\claimskip\noindent{\bf#1.}
{\it\ignorespaces#2}\vskip\claimskip\noindent}

\font\ninerm=cmr9
\font\ninei=cmmi9
\font\ninesy=cmsy9
\font\ninebf=cmbx9

\font\nineit=cmti9
\font\ninesl=cmsl9

\font\eightrm=cmr8
\font\eighti=cmmi8
\font\eightsy=cmsy8
\font\eightbf=cmbx8

\font\eightit=cmti8
\font\eightsl=cmsl8
\font\sixrm=cmr6
\font\sixi=cmmi6
\font\sixsy=cmsy6
\font\sixbf=cmbx6
\font\fiverm=cmr5
\font\fivei=cmmi5
\font\fivesy=cmsy5
\font\fivebf=cmbx5

\newfam\itfam
\newfam\slfam
\newfam\bffam
\def\eightpoint{\def\rm{\fam0\eightrm}%
\textfont0=\eightrm \scriptfont0=\sixrm \scriptscriptfont0=\fiverm
 \textfont1=\eighti \scriptfont1=\sixi \scriptscriptfont1=\fivei
 \textfont2=\eightsy \scriptfont2=\sixsy \scriptscriptfont2=\fivesy
 \def\it{\fam\itfam\eightit}%
 \textfont\itfam=\eightit
 \def\sl{\fam\slfam\eightsl}%
 \textfont\slfam=\eightsl
 \def\bf{\fam\bffam\eightbf}%
 \textfont\bffam=\eightbf \scriptfont\bffam=\sixbf
 \scriptscriptfont\bffam=\fivebf
 \normalbaselineskip=9pt
 \setbox\strutbox=\hbox{\vrule height7pt depth2pt width0pt}%
 \normalbaselines\rm}

\def\ninepoint{\def\rm{\fam0\ninerm}%
\textfont0=\ninerm \scriptfont0=\sixrm \scriptscriptfont0=\fiverm
 \textfont1=\ninei \scriptfont1=\sixi \scriptscriptfont1=\fivei
 \textfont2=\ninesy \scriptfont2=\sixsy \scriptscriptfont2=\fivesy
 \def\it{\fam\itfam\nineit}%
 \textfont\itfam=\nineit
 \def\sl{\fam\slfam\ninesl}%
 \textfont\slfam=\ninesl
 \def\bf{\fam\bffam\ninebf}%
 \textfont\bffam=\ninebf \scriptfont\bffam=\sixbf
 \scriptscriptfont\bffam=\fivebf
 \normalbaselineskip=11pt
 \setbox\strutbox=\hbox{\vrule height7pt depth2pt width0pt}%
 \normalbaselines\rm}

\def\plainsection#1{\par\vskip .5cm\penalty -100 
\vbox{\noindent{\sc #1}
\vskip 5pt}
\penalty 500}

\def\Bibitem#1&#2&#3&#4&%
{\hangindent=1.66cm\hangafter=1
\noindent\rlap{\hbox{\rm #1}}\kern1.66cm{\rm #2}{\it #3}{\rm #4.}
\vskip5pt}

\def\square{{\hfill \hbox{
\vrule height 1.453ex  width 0.093ex  depth 0ex
\vrule height 1.5ex  width 1.3ex  depth -1.407ex\kern-0.1ex
\vrule height 1.453ex  width 0.093ex  depth 0ex\kern-1.35ex
\vrule height 0.093ex  width 1.3ex  depth 0ex}}}
\def\bigsquare{{\kern-0.3ex\hbox{
\vrule height 1.7ex  width 0.093ex  depth 0ex\kern-0.093ex
\vrule height 1.8ex  width 1.7ex  depth -1.707ex\kern-0.093ex
\vrule height 1.7ex  width 0.093ex  depth 0ex\kern-1.65ex
\vrule height 0.093ex  width 1.6ex  depth 0ex}\kern0.3ex}}
\def\qed{\phantom{~}$\square$\medskip}
\def\smallskip{\vskip 3pt}
\def\medskip{\vskip 5pt}

\title{A class of strictly pseudoconvex domains with non-pluripolar core}

\author{Zbigniew Slodkowski}
\date{July 2021}

\begin{document}


\begin{abstract}
We construct a class of strictly pseudoconvex domains in $\bC^d$ whose core has non-empty interior. Consequently these cores are not pluripolar. This answers a question posed by Harz, Shcherbina and Tomassini.
\end{abstract}

\maketitle 

\plainsection{0. Introduction} 

The {\it core} of a complex manifold $M$ (denoted ${\bf c} (M)$ )  is the largest set on which every bounded and continuous plurisubharmonic function on $M$ fails to be strongly plurisubharmonic. This notion was introduced and studied by Harz, Shcherbina and Tomassini in [HST17, 20, 21], originally for unbounded, strictly pseudoconvex domains in $\bC^d$, and then also for arbitrary complex manifolds, cf. [PoS19]. In this generality the concept of the core blends somewhat with the notion of minimal kernel introduced earlier in [SlT04], cf. [MST18], defined in terms of not-necessarily bounded plurisubharmonic functions of a fixed regularity class. Similarly, for different regularity classes, corresponding notions of the core can be defined (potentially different, but demonstrating examples are not known), cf. [HST17], [Slo19]. The smallest of them is the one defined with respect to all bounded continuous plurisubharmonic functions ({\it the continuous core}), and we restrict our discussion here to this case. Also, we consider here only the core of strictly pseudoconvex domains in $\bC^d$ (necessarily unbounded for the core to be nonempty).

In [HST17], the authors ask a series of questions about structure of the core. Answering one of them (Question 2, Section 4) we prove here that there is a strictly pseudoconvex domain in $\bC^2$ whose core is not pluripolar, and moreover, has nonempty interior.

The general idea behind our examples is to consider a closed subset $Z$ of $\bC^2$, with nonempty interior, which is contained in the limit (in a certain sense) of a sequence $(M_j)$ of parabolic Riemann surfaces, and such that the union
$$
Z\cup \bigcup_j M_j \leqno (0.1)
$$
is contained in some strictly pseudoconvex domain $\Omega$. Since every bounded and continuous plurisubharmonic function on $\Omega$ must be constant on every $(M_j)$, their union, and its closure, are contained in the core of $\Omega$. Consequently, the core contains $Z$ (which has nonempty interior). In order that such $\Omega$ exists, $Z$ must satisfy some conditions: we assume that there is a continuous plurisubharmonic function $\lambda$ on $\bC^2$, bounded on $Z$, and such that its all sublevel sets are hyperconvex. (The author hopes that this class of functions might be useful in some other problems).

The general results outlined above (with necessary definitions) are stated in Section 1, and proven in Section 2. In Section 3, as a specific example of set $Z$ we use the "tubular" set bounded by the Levi-flat hypersurface constructed by Brunella, with the help of a Fatou-Bieberbach mapping (in response to another question posed in [HST12]), and the special plurisubharmonic function $\lambda$ mentioned above is constructed from the coordinate functions of that F-B mapping. In Section 4, we give some variations of the basic example, in particular (Remark 4.5) relating it to a different notion of the core due to Gallgher et al [GHH17]. Finally, in the Appendix, we sketch proofs of some undoubtedely known facts, for which we could not find a reference.

\claim 0.1. Some terminology and notation|\rm We list now terms not defined in the text, specifying, in cases of ambiguity, the meaning we assume here.  
\endclaim
  
We abbreviate {\it plurisubharmonic} to PSH; thus $PSH(W)$ and $C(W)$ stand, respectively for classes of all PSH or continuous functions on the domain $W$. 

For us {\it strictly} PSH functions are at least $C^2$ -smooth; a PSH function $\phi$ (that does not have to be continuous) is {\it strongly} PSH in a domain 
$W \subset \bC^d$, if for every $\rho\in PSH(W)$ with compact support in $W$, the function $\phi + t\rho \in PSH(W)$, for $|t|$ small enough. 

Similarly, a {\it strictly} pseudoconvex domain must have at least 
$C^2$ -smooth boundary.

A domain (possibly unbounded) is called {\it hyperconvex} if it has a {\it bounded continuous} PSH exhaustion function, [KeR81].

A complex variety in $\bC^d$ is called {\it parabolic} if every upper-bounded PSH function on it must be constant, cf. [AyS14].

$D$ will denote the unit disc.

\claim  Acknowledgement|\rm I am grateful to J.-P. Demailly for providing references regarding parabolic varieties.  
\endclaim 

\plainsection{1. Surrounding a sequence of varieties by a strictly pseudoconvex hypersurface}

In this section we present the basic theorems on which our construction of the example is based. The proofs will be given in Section 2.

\claim 1.1. Theorem|Let $Z$ be a closed subset of $\bC^d$, and $M_j, j\in \bN$, a sequence of irreducible varieties of codimension one, each defined by an entire holomorphic function $h_j$ on $\bC^d$. Assume that for every $p\in \bC^d \setminus Z$ a sufficiently small ball $B(p,\epsilon), \epsilon >0$, intersects at most finitely many $M_j$'s. Assume further that there is a nonnegative continuous PSH function $\lambda$ on $\bC^d$ which is bounded on $Z$ and such that its all (nonempty) sublevel sets $\{z\in\bC^d : \lambda (z) < C\}$, $C > 0$, are hyperconvex domains.  Assume finally that there is a positive PSH function $\mu$ on $\bC^d$ which is bounded on $Z$ and strictly PSH on $\bC^d \setminus Z$. Then there is a PSH function $\psi$ on $\bC^d$ such that

\emph{(i)}  $\psi|M_j \equiv -\infty$, for $j\in\bN$;

\emph{(ii)} $\psi$ is bounded from the above on $Z$;

\emph{(iii)}  $\psi$ is $C^\infty$ smooth and strictly PSH on $\bC^d\setminus (Z\cup \bigcup_j M_j)$.

\endclaim

The assumption about existence of additional function $\mu$ could be dropped by requiring that the function $\lambda$ is, in addition, strictly PSH. Whether this would result in {\it actual} loss of generality is not clear, but it would suffice for our specific application. The proof of the theorem is based on the next lemma.

\claim 1.2. Lemma| Let $W \subset\bC^d$ be a (possibly unbounded) hyperconvex domain, $h : W \to \bC$ a holomorphic function with nonempty zero set $M$, and $E\subset M$ a compact set. Then for every $\delta > 0$, and for every compact subset $ L$ of $ \overline{W}\setminus M$ there is a continuous function $v : \overline{W} \setminus M \to [-1, 0)$, such that

\emph{(i)}  $v|E \equiv -1$;

\emph{(ii)} $v|L > -\delta$;

\emph{(iii)}  $\sup v | \overline{W} < 0$.

\endclaim

The above theorem gives a way to include some sequences of varieties in a strictly pseudoconvex domain (of the form $\{\psi < C\}$). These varieties have to {\it cluster} on the set $Z$ in the following sense.

\claim 1.3. Definition| Let $Z$ and $X_j$, for $j\in\bN$, be nonempty closed subsets of an Euclidean space. We say that $Z$ is the cluster of the sequence $(X_j)$, if the following two conditions hold:

\emph{(i)}  whenever $(p_k)$ is a convergent sequence of points, such that $p_k\in X_{j_k}$, where $j_k < j_{k+1}, k\in \bN$, then $\lim p_k \in Z$;

\emph{(ii)} every $p\in Z$ is the limit of a sequence $(p_k)$ like in (i).

\endclaim
In particular, for $p$ not in $Z$ a sufficiently small ball $B(p,\epsilon), \epsilon >0$, intersects at most finitely many $V_j$'s, the condition 
required in Theorem 1.1. In our example, we will have to apply Theorem 1.1 when $M_j$'s are parabolic varieties and $Z$ is of the form
$$
Z:= \{ (z,w)\in \bC^2 : |g(z,w)| \leq 1 \},
\leqno(1.4)
$$
where $g(z,w)$ is an entire holomorphic function. To apply Theorem 1.1 we need to construct a sequence $(M_j)$ of {\it parabolic} varieties that clusters on $Z$. We will do it for a more general class of sets $Z$. 

\claim 1.5. Lemma| Let $Z \subset \bC^d$ be the cluster of a sequence $(V_j)$ of irreducible varieties of codimension one. Then $Z$ is the cluster of some sequence $(M_j)$ of of irreducible parabolic varieties of codimension one.

\endclaim
Since sets (1.4) are foliated by varieties $\{g = const\}$, the assumption of the lemma are satisfied. The following statement summarizes the method of construction of our class of examples.

\claim 1.6. Corollary| Let $Z$ be the cluster of a sequence of irreducible parabolic varieties of codimension one in $\bC^d$ satisfying further assumptions of Theorem 1.1. Then there exist strictly pseudoconvex domains $\Omega$ containing set (0.1), and the core of every such domain contains the set (0.1).

\endclaim{}
\proof Let $\psi$ be the PSH function constructed in Theorem 1.1. In particular  $\psi$ is $C^\infty$ smooth and strictly PSH on $\bC^d \setminus (Z\cup \bigcup_j M_j)$. Furthermore, $C_0: = \sup \psi | Z$ is finite, and so $\psi$ is $\bC^\infty$ smooth on the domain $\{ z \in \bC^d : \psi(z) > C_0\}$. By Sard theorem, almost every $C > C_0$ is a regular value of $\psi$. For such $C$ the (unbounded) domain $\Omega := \{z \in \bC^d : \psi (z) < C  \}$ is strictly pseudoconvex.
Let $u$ be any bounded continuous PSH function on $\Omega$. Then $u|M_j$ is a constant function for every j (as $M_j$ is parabolic). Thus there is no $p \in M_j$ and $r > 0$ such that $u | B(p, r)$ is strongly PSH. Since core 
${\bf c} (\Omega)$ is a closed set (in $\bC^d$, when $\Omega$ is strictly pseudoconvex), and 
$Z \subset \overline{\bigcup_j M_j}$, we get $Z \subset {\bf c} (\Omega)$. \qed

\plainsection{2. Proofs of Theorem 1.1 and Lemma 1.5.} 

Crucial to what follows Theorem 1.1 belongs to the long sequence of results descending from Josefson Theorem, and, in its general scheme, it's proof is similar to that of Bedford and Taylor [BeT82]. We follow here the transparent presentation in Kolodziej [Kol05] (Theorem 1.23) which starts from a convergent series of constant multiples of the function 
$\log (1 +|z|^2)$, and then these multiples are modified on progressively larger balls, so that the sum of the resulting series is equal to $-\infty$ on the given union of pluripolar sets. In our construction  
$\log (1 +|z|^2)$ is replaced by $\lambda (z)$, and the sequence of Euclidean balls by a sequence of progressively larger sublevel sets of $\lambda$. The semi-local modification step is more complicated here because we need to assure that our series is locally uniformly convergent on a large subdomain of $\bC^d$.

\claim 2.1 Local construction|
\endclaim

\noindent{\it Proof of Lemma {\rm 1.2.}} As $W$ is hyperconvex, it has a bounded exhaustion function $\rho ^* \in PSH (W) \cap C(\overline{W})$ with $\sup \rho ^* = 0$. For some $c > 0$, $c\rho ^* |E \leq -1$. Let 
$$
\rho (z): = \max (c\rho ^*(z) - \delta /4, -1), \hbox{ , for } z \in W, \hbox { , and }\rho (\zeta): = -\delta /4 , 
\hbox { , for } \zeta \in bW. 
$$
 
Function $\rho$ has the following properties:

{\plainitemindent=14mm
\plainitem{(2.1.1)} $\rho\in PSH(W) \cap C(\overline{W})\;$;

\plainitem{(2.1.2)} $\rho : \overline{W} \to [-1, -\delta / 4]\;$ and $\rho |E = -1\;$;

\plainitem{(2.1.3)} if $(z_s)$ is a sequence in $W$ , and $\lim |z_s| = +\infty$, then $\lim\rho (z_s) = -\delta /4\;$;

\plainitem{(2.1.4)} for every $t\in [-1, -\delta /4)$ the subset $\{z\in W : \rho (z) \leq t \}$ is compact.

\noindent Let $H_1: =\{z \in W : \rho(z< -\delta /2\}$ and $H_2: = \{z \in W : \rho (z) < -3\delta /4\} $. Then $H_1$ and $H_2$ are relatively compact open sets, such that $\overline{H}_2 \subset H_1$ and $\overline{H}_1\subset W$. $\overline{H}_1$ being compact, let $m: = \max \{|h(z)| : z \in \overline{H}_1\} $.
Let $L^*:= L\cap\overline{H}_1$. 
As $\log (|h(z)|/m)$ is a continuous PSH function on $W\setminus M$ with maximum value $0$ on $\overline{H}_1$ and finite minimum value on the compact $L^*$, there is a positive constant $\kappa$ such that $\kappa (\log |h(z)|/m) > -\delta /4$, for $z \in L^*$ Thus if
$$
\phi (z): = \kappa (\log |h(z)|/m) - 3\delta /4, \hbox{ for } z \in H_1,
$$

\noindent then $\phi \in PSH(H_1)\cap C(H_1\setminus M)$, and satisfies $ \phi(z) > -\delta $ on $ L^* $ and $\phi (z) < - 3\delta /4$ on $H_1$. 
Finally define $v$ on $\overline{W}$ by
$$
v(z): = \begin{cases}
        \rho(z), & \text { if  $z\in \overline{W} \setminus H_2$}\\
        \max(\rho(z), \phi(z)), & \text { if $z\in H_2 $}.
        \end{cases}
$$

\noindent Since $\rho(z)\geq -3\delta/4 $ on $\overline{W}\setminus H_2 $, and $\phi(z) < -3\delta/4$ on $H_1$, so $\max(\rho,\phi) = \rho $ on $H_1\setminus \overline{H}_2$. Thus $v\in PSH(W)$, and
$$
v\leq -\delta/4  \hbox{ on } \overline{W}. \leqno (2.1.5)
$$
\noindent As $\phi \in C(H_1 \setminus M)$, and $\rho > \phi$ on $H_1\setminus \overline{H}_2$,
we get $v \in C(\overline{W})$, as well as that 
$$
v|E = -1\hbox{ and } v\geq -1 \hbox{ on } \overline{W} \leqno (2.1.6)
$$
\noindent (using also (2.1.2)). As $v| L^* > -\delta$, and since $L\subset L^*\cup (\overline{W}\setminus H_1)$ and $\rho\geq -\delta/2$ on $\overline{W}\setminus H_1)$, we obtain 
$$
v|L > -\delta.\leqno (2.1.7)
$$
Relations (2.1.5), (2.1.6), (2.1.7) are the properties (i),(ii),(iii) required in Lemma 1.2.
\qed

\claim 2.2 Special plurisubharmonic function|
\endclaim

\noindent {\it Proof of Theorem {\rm 1.1.}} Choose $k_0 \in \bN$ such that $\sup \lambda |Z < 2^{k_0}$ , and let $W_k: = \{z \in \bC^d : \lambda (z) < 2^k\}$, for $k \geq k_0$.
By assumptions of the theorem, $W_k$'s are hyperconvex, and so we can apply Lemma 1.2 to produce a function $v_k \in PSH(W_k)\cap C(\overline{W}_k)$, with special properties. We first need to specify, for each $W_k$, the sets corresponding to $L$ and $E$ in the lemma.

Due to the assumption of local finiteness in $\bC^d \setminus Z$ of the configuration of $M_j$'s, the union 
$Z \cup M^*$, where $M^*:= \bigcup_j M_j$, is closed in $\bC^d$.
Choose a convergent to $0$, strictly decreasing sequence $(\epsilon_k) \subset (0,1)$, and define open sets 
$U_k: =\{z\in \bC^d : dist (z, Z\cup M^*) < \epsilon_k \}$, for $k\geq k_0$. Clearly
$$
\overline{U}_{k+1} \subset U_k \hbox {, for } k\geq k_0 , \hbox{ , and } \bigcap_{k\geq k_0} U_k = Z\cup M^* .\leqno(2.2.1)
$$
Consider also a sequence of concentric open balls $B_k$, for $k > k_0$, such that
$$
\overline{B}_k \subset B_{k+1} \hbox{ , with } B_k\not\subset\overline{W}_k \hbox{, and } \bigcup_{k>k_0} B_k = \bC^d .\leqno(2.2.2) 
$$
\noindent Let finally $L_k: =\overline{W}_k\cap\overline{B}_k\setminus U_k $ , for $k > k_0$ . We can assume WLOG (diminishing $\epsilon _k$'s , if needed) that all $L_k$'s are nonempty.

To define sets $E$, choose, for each $j$ a compact subset $E_j$ of $M_j$, with non-empty interior relative to $M_j$. In the same way as in [BeT82], [Kol05] we want that the same $E_j$ is used in infinitely many $W_k$'s, and that when it is used as $E$, for a $W_k$, it is fully contained in it, setting it as empty , if this is not possible. Cutting short the combinatorial details, we claim that there is a function $j : \{ k > k_0\} \to \{0\}\cup\bN$, such that each $j \in \bN$ appears as $j(k)$ for infinitely many $k$'s, and 
$$
E_0 =\emptyset \hbox{ , and } E_{j(k)}\subset W_{k-1} \hbox { , for } k > k_0 .
$$

We apply now Lemma 1.2 with $W:= W_k$ (which is hyperconvex), $h:= h_k|W_k$ , $M:= M_{j(k)}\cap W_k$ , $E:= E_{j(k)}$ , and $L:= L_k$, for $k > k_0$. By (2.2.1), (2.2.2), we have $L\subset \overline{W}\setminus M$. If $j(k)=0$ we let $v_k = -2^{-k-1}$ on $\overline{W}_k$. Otherwise applying the lemma with specified data and $\delta := -2^{-k}$ for $k > k_0$, we obtain a function 
$v_k\in C(\overline{W}_k)\cap PSH(W_k)$, such that
$$
v_k|E_{j(k)} \equiv -1 \hbox{ , and } v_k\geq-1 \hbox{ , on } \overline{W}_k \hbox{ ;}
$$
$$
\sup v_k| \overline{W}_k < 0 ; \hbox{ , and }  v_k| L_k > -2^{-k} .\leqno (2.2.3)
$$

We will construct now an entire PSH function $\phi$ by {\it gluing} $v_k$'s  to $\lambda$. For $k > k_0$ let
$$
\phi_k (z): = \begin{cases}
              \max (v_k(z),  2^{-k+1}\lambda(z)-2), & \text { if  $z\in W_k$} \\
              2^{-k+1}\lambda(z)-2,  & \text { if $z\in \bC^d \setminus W_k$},
              \end{cases}
$$
and then, for $z\in\bC^d$, let
$$
\phi(z):= \sum_{k > k_0} \phi_k(z).
$$
Observe first that $2^{-k+1}\lambda(z)-2 = 0$  on $bW_k$, while $\sup v_k| \overline{W}_k < 0$, and so
$\phi_k \in PSH(\bC^d)\cap C(\bC^d)$.
Now, for any fixed $k_1 > k_0$, we have $\phi_k |W_{k_1} \leq 0$, for all $k\geq k_1$, hence $\sum_{k\geq k_1} \phi_k$ is PSH on $W_{k_1}$, and so is $\phi|W_{k_1}$ (for every $k_1$). Thus $\phi\in PSH(\bC^d)$, provided it is not identically equal to $-\infty$, which we address now.

By (2.2.3),  $|\phi_k (z)|\leq 2^{-k}$, for $z\in L_k$, and $k >k_0$. Since $(L_k)$ is an increasing sequence of sets (by (2.2.1),(2.2.2)), the series of continuous functions $\sum_{k\geq k_1} \phi_k$ is uniformly convergent on $L_{k_1}$ for every $k_1 > k_0$. Hence 
$\phi$ is a continuous function on each 
$L_{k_1}$.
Consider now an arbitrary point $p \not\in Z\cup M^*$. By (2.2.1), (2.2.2), there is a ball $B(p,r)$, with $r>0$ contained in some $L_k$. Thus $\phi$ is continuous on that ball, and consequently on 
$\bC^d\setminus (Z\cup M^*)$.
We will show now that
$$
\phi|M^* \equiv -\infty.
$$
Consider any set $E_j$. Then $E_j$ appears infinitely many times as $E_{j(k_s)}$, where  $k_1 < k_2 < ...$.
Let $p \in E_j$. For any $k = k_s$, as $p$ belongs to $W_{k-1}$, we have $v_k(p) = -1$, and $2^{-k+1}\lambda(p) -2 \leq -1$, and so $\phi_{k_s} (p) = -1$, for all $s\in \bN$. By similar reasons, $\phi_k(p)\leq 0$ for all $k\geq k_1$. Thus $\phi(p) = -\infty$ for $p \in E_j$. Since, (for every $j$), $E_j$ has nonempty (relative) interior in $M_j$ (which is irreducible),
$\phi |M^*\equiv -\infty$.

Note that, as $v_k < 0$ on $W_k$, and $\lambda$ is nonnegative, $\phi_k (z)\leq 2^{-k+1}\lambda (z)$, and so 
$$
\phi(z)\leq 2\lambda(z) \hbox { , on  }  \bC^d. \leqno (2.2.4)
$$

To finish, we regularize $\phi$. Observe $\phi +\mu$ is PSH on $\bC^d$, and continuous and {\it strongly} PSH on 
$\bC^d\setminus (Z\cup M^*)$. We use the {\it tangential} approximation of Richberg [Ric68], [Smi86]. Thus, let $\epsilon (z)$ be a continuous strictly positive function on $H:=\bC^d\setminus (Z\cup M^*)$ with $\lim_{z \to bH}\epsilon(z) =0$, and $\epsilon(z) < \lambda(z)$, on $H$. Then there is $\psi\in PSH(\bC^d)$, strictly  PSH on H, and such that
$$
\psi(z)=\phi(z) + \mu (z) \hbox{ on } Z\cup M^*;
$$
$$
(\phi + \mu)(z)\leq \psi (z) \leq (\phi+\mu+\epsilon)(z)  \hbox{ , for } z\in H.
$$
Observe that by (2.2.4) $\psi(z)\leq  3\lambda(z) + \mu (z)$  on  $Z\cup M^*$ and so $\psi$ is bounded from the above on $Z$ and is identically $-\infty$ on $\bigcup _j M_j$. Thus conditions (i) - (iii) of Theorem 1.1 hold. 
\qed

\claim 2.3 Approximations by parabolic varieties|
\endclaim

\noindent {\it Proof of Lemma {\rm 1.5.}} Choose a sequence $(F_k)$ of compact subsets of $\bC^d$, such that
$$
F_k\subset Int F_{k+1} \hbox { , and  } \bigcup_k F_k = \bC^d\setminus Z. \leqno(2.3.1)
$$
Observe first that for every $k$ there can be at most finitely many varieties $V_j$ intersecting $F_k$. (Otherwise one could find a convergent sequence $(p_n)$ with $p_n \in F_k \cap V_{j_n}$, where $j_n < j_{n+1}$. As $\lim p_n \in F_k$, it cannot belong to $Z$, contradicting Definition 1.3(i).)

Choose a countable dense subset $A$ of $Z$ and put its points in a sequence $(a_j)$ in such a way that every $a \in A$ appears as $a_j$ for infinitely many $j$'s. We will construct now a sequence of parabolic varieties $M_j \subset \bC^d$, and a sequence of points $(p_j)$ satisfying, for all $j$,
$$
M_j \cap F_j = \emptyset  ;\leqno(2.3.2)
$$
$$
p_j \in M_j; \leqno(2.3.3)
$$
$$
|p_j - a_j| < 1/j. \leqno (2.3.4)
$$
Fix an arbitrary $j \in \bN$, and choose $r \in (0, 1/j)$ such that $F_j \cap \overline{B(a_j, r)} = \emptyset $.
Using the observation that given $F_j$ can intersect at most finitely many varieties $V_k$, and recalling they satisfy condition (ii), Definition 1.3, we conclude that there is a variety $V_{k(j)}$ and a point $p_j$ such that
$$
p_j \in V_{k(j)} \cap B(a_j, r) \hbox{ , and  } V_{k(j)} \cap F_j = \emptyset.
$$
Now $V_{k(j)}$ is defined by some entire holomorphic function: $V_{k(j)} = \{ z \in \bC^d : g(z) = 0 \}$. As $V_{k(j)}$ is disjoint from $F_j$,we have
$$
0 < m:= \min \{|g(z)| : z \in F_j \}.\leqno (2.3.5)
$$
We can approximate the entire function $g$  uniformly on the compact set $F_j \cup \{p_j\}$ by a polynomial $P(z)$ so that
$$
|g(z) - P(z)| < m/4 \hbox{ , for } z\in F_j\cup \{p_j\}.
$$ 
As $g(p_j) = 0$, we have $|P(p_j)| < m/4$, and so 
$$
|g(z) - (P(z)-P(p_j))| < m/2 \hbox{ , for } z\in F_j.
$$
By this and (2.3.5),  $|P(z)-P(p_j)| > m/2 $ on $F_j$. Setting as $M_j$ the irreducible component (containing $p_j$)  of the variety  $\{P(z) -P(p_j) = 0\}$, we can see that so constructed $M_j$ and $p_j$ satisfy conditions (2.3.2 - 4). 
In particular, (2.3.3) and (2.3.4) imply, together with the fact that points of $a_j$ are repeated infinitely many times each, and are dense in $Z$, that condition (ii) of Definition 1.3 holds. As well known, $M_j$'s must be parabolic (For example, this foloows from Theoreme 7.5 in [Dem85].)

To verify that condition (i) of that definition holds, consider a convergent sequence $(p_k)$ such that $p_k \in M_{j(k)}$, where $j(k) < j(k+1)$, with $p:=\lim p_k $. Suppose, to the contrary, that $p \not\in Z$. By (2.3.1), $p \in Int F_{j^*}$, for some $j^*$ and so infinitely many of $M_{j(k)}$'s intersect $F_{j^*}$. This is a contradiction, because by (2.3.2) and  (2.3.1), only finitely many of $M_j$'s (with $j < j^*$) can intersect $F_{j^*}$. 
\qed

\plainsection{3. Application of Fatou-Bieberbach domains}

We will conclude now, just in $\bC^2$, the construction of a strictly pseudoconvex domain $\Omega$ with {\it thick} core, by giving an explicit example of set $Z$, to which Corollary 1.6 applies, and of a function $\lambda$, as in Theorem 1.1 (on which Corollary 1.6 is based). In fact, the boundary of our set $Z$ is the Levi flat hypersurface in $\bC^2$ constructed by Brunella in response to another question of Harz et all, cf. [HST12]. Brunella's construction was based on the example of a Fatou-Bieberbach domain due to Globevnik [Glo98], from which we distill only the features we need.

\claim 3.1. Theorem [Glo98]| There is a domain $U$ in $\bC^2$ and a biholomorphic map $\Phi = (f,g)$ of $\bC^2$ onto $U$, such that

\emph{(i)}  $U \subset \{(\xi,\zeta) \in \bC^2: |\xi| < \max (C,|\zeta|) \} $, for some constant $C$ ;

\emph{(ii)} the projection of $U$ on second coordinate fills $\bC$.

\emph{(iii)} $U$ has $C^1$-smooth boundary.

\endclaim

\claim 3.2. Proposition|With $\Phi = (f,g)$ as in Theorem 3.1), let $Z:= g^{-1} (\overline{D})$, and let  
$$
\lambda(z,w):= \log (1 + |\Phi(z,w)|^2) \hbox{ , for } (z,w) \in \bC^2 .
$$
Then $Z$, and functions $\lambda$, and $\mu := \lambda$, satisfy assumptions of Lemma1.5 and Theorem 1.1.
\endclaim

\proof To show that $Z$ satisfies conditions of Lemma 1.5, we choose a dense sequence $(b_n)$ of points in disc $D$, and let $V_n$
be varieties $\{g = b_j\}$. To show that $Z$ is the cluster of the sequence $(V_j)$, note first that the condition (i) of Definition 1.3 is satisfied trivially, because all $V_j$'s are contained in $Z$. To verify (ii), consider any $(z_0, w_0) \in Z$. Then $(a,b) := \Phi(z_0, w_0)$ belongs to  $U \cap (\bC\times\overline{D})$. Clearly, for some increasing sequence of indices $(j(n))$, we have 
$(a,b) = \lim (a,b_{j(n)})$ with $(a,b_{j(n)})\in U\cap (\bC \times\overline{D})$. Letting now $p_n :=\Phi^{-1}(a,b_{j(n))}$, we obtain $\lim p_n = (z_0,w_0)$, with $p_n\in V_{j(n)}$, which confirms condition (ii) of Definition 1.5.

By Theorem 3.1 (i), $\Phi(Z) = U \cap \{(\xi,\zeta) : |\zeta|\leq 1\} \subset C\overline{D}\times\overline{D}$, and so $\lambda$ is uniformly bounded on $Z$ (by $\log (2 + C^2)$).

It is clear that $\lambda$ is $C^{\infty}$-smooth and non-negative everywhere. It is strictly PSH as composition of the strictly PSH function $\log(1 + |\xi|^2 + |\zeta|^2)$ with the biholomorphic map $\Phi$. 

It remains to show that the sublevel sets $\{\lambda < C\}$ , $C > 0$ , are hyperconvex. 
We use the notion of {\it local} hyperconvexity. As in [KeR81], p. 172] we say that a domain $H$ is locally hyperconvex at 
$a\in bH$, if there is $r > 0$ and $\phi \in PSH(H\cap B(a,r))\cap C(H\cap B(a,r))$, such that $\lim \phi (z) = 0$, as  $z \to bH$, for $z$ in $H \cap B(a,r)$. Such $\phi$ is called a local barrier function at $a$ (for $H$). By [KeR81], Lemma 2, if $H$ is a (possibly unbounded) domain with $C^1$-smooth boundary, then $H$ is locally hyperconvex at every boundary point.
Now if a domain $H = H_1 \cap H_2$, and point $a \in bH$, then $a \in bH_1$ and $a \in bH_2$. 
If, for $i=1, 2$, domain $H_i$ is locally hyperconvex at $a$, with barrier function $\phi_i$, then $H_1 \cap H_2 = H$ is locally hyperconvex at $a$, with the barrier function $\phi := \min (\phi_1,\phi_2)$. Applying these observations to $H_1:= U$, which is pseudoconvex with $C^1$ boundary (by Theorem 3.1), and to the obviously locally hyperconvex sphere 
$H_2:=\{\log(1 + |\xi|^2 + |\zeta|^2) < C\}$, we obtain that their intersection is locally hyperconvex. But since this intersection is relatively compact in $\bC^2$, it must be hyperconvex, by [KeR81], Proposition 1.1]. The image of this intersection by the biholomorphic map $\Phi ^{-1}$ is 
$\{\lambda < C \}$, and so the latter set is hyperconvex as well (hyperconvexity being biholomorphically invariant). 

It follows that $\lambda$, and $\mu = \lambda$, have all the required properties.
\qed

\claim 3.3. Conclusion: the example|\rm By Proposition 3.2 and Lemma 1.5, the set $Z$ is the cluster of some sequence $(M_j)$ of parabolic varieties, which, together with functions $\lambda$, and $\mu = \lambda$, can be used in Theorem 1.1, to obtain a special plurisubharmonic function $\psi$, which, by Corollary 1.6, yields a strictly pseudoconvex domain $\Omega$ whose core contains $Z$, a set with nonempty interior, and so nonpluripolar.
\endclaim

\claim 3.4. Remark|\rm Function $\lambda$ can be chosen as in Proposition 3.2, or just as $\lambda(z,w):= |\Phi(z,w)|^2 $. The same two choices can be done for $\mu$, and independently of those for $\lambda$. They would all produce an example like in 3.3, but if we want  $\Omega$ to have extra properties, the choice of $\mu$ (or $\lambda$) could matter. See Remark 4.5 for such situation.
\endclaim

\plainsection{4. Related examples and observations}

\claim 4.1. Version of the problem for components of the core|\rm Considering the possibility that the core might be, in general, not pluripolar, Nikolai Shcherbina asked (during a conversation in the Fall 2017), whether, at least, its {\it components} are pluripolar. However, the answer to this question is also negative.
\endclaim

\noindent Specifically, recall that for a point $p\in {\bf c}(\Omega)$, the component of the core containing $p$, is, in the notation of [PoS19], the set $A^b_e (p)$ of all points $z\in\Omega$, such that $u(z) = u(p)$ for all  bounded $u \in PSH(\Omega)\cap C(\Omega)$. Components were introduced in [HST17], cf. also [PoS19] and [Slo19] for further information. To construct an example of a {\it thick} component consider the set $Z$ and the sequence $(M_j)$ of parabolic varieties, with which example of $\Omega$ was build in Section 3. By Proposition A.1, for each $M_j$, there are at most finally many exceptional complex directions in $\bC^2$, such that some complex line parallel to one of them does not intersect $M_j$. Consequently, there is a complex line that intersects all varieties $M_j$, where $j \in \bN$; denote this line by $M_0$. Now, if we reapply construction of Sections 1 - 3 to $Z$ and $\tilde {M}:= \bigcup_{j\geq 0} M_j$, we obtain a strictly pseudoconvex domain $\tilde{\Omega}$ containing 
$Z \cup \tilde{M}$. Fix a point $a \in M_0$, and let 
$u \in PSH(\tilde{\Omega})\cap C(\tilde{\Omega})$ be a bounded function. Then $u$ is constant on each $M_j$, and so also on $\tilde{M}$, as it is connected. Since $Z$ is contained in the closure of  $\tilde{M}$, $u|Z = u(a)$.
Thus $Z \subset A^{cb}_e (a)$, i.e. the component containing $a$ has nonempty interior.

\claim 4.2. A finer decomposition of the core|\rm  Poletsky and Shcherbina [PoS19] considered also another relation in $\Omega$, where $p, q$ are equivalent, if they are not distinguished by any bounded function in $PSH(\Omega)$. The equivalence class containing $p$ is denoted $A^b_e(p)$. For lack of a better term we call these classes {\it strata}. They gave ([PoS19, Section 5] an example showing that, anlike components, a stratum does not have to be closed, and a stratum contained in the core can be a single point. Their domain $\Omega \subset \bC^2$ is {\it not} strictly pseudoconvex (and its core is not closed in $\bC^2$); we observe now that the same phenomena occur in our strictly pseudoconvex $\Omega$ constructed in Section 3.
\endclaim

\claim 4.3. Observation|\rm Let $Z$, $(M_j)$, $\lambda$ be like in Section 3, with $\mu$ assumed strictly PSH on $\bC^2$). and $\Omega$ be like in Section 3. Assume further that $M^* = \bigcup_j M_j$ is connected. Then

\emph{(i)}  the stratum $S$ containing $M^*$ is pluripolar and not closed (neither in $\Omega$, nor in $\bC^2$);

\emph{(ii)} every point in $\Omega \setminus \psi^{-1}(-\infty)$ is a stratum in itself.
\endclaim

\proof \rm (i) With $\psi$ as in Theorem 1, $\psi |M^* = -\infty$, and so 
$S$, contained in $\psi^{-1}(-\infty)$, is pluripolar. Since 
$$
M^*\subset S \subset \psi^{-1}(-\infty) \subset M^*\cup Z  \hbox{ , and }   \overline {M^*} = Z\cup M^*), 
$$
$S$ is dense in $M^*\cup Z$, but $S\neq M^*\cup Z$, because $Z$ is not pluripolar. 
Hence $S$ cannot be closed.

\rm (ii) Fix any point $a\in \Omega \setminus \psi^{-1}(-\infty)$, and consider an arbitrary point $b$ in $\Omega\setminus {a}$. If $\psi (b) = -\infty$ then $\psi (b) \neq \psi(a)$, and we are done. Otherwise both $\psi (b) , \psi(a)$ are finite, which, by construction of $\psi$ (in the proof of Theorem 1.1), implies that $\alpha (a), \alpha (b)$ are both finite, where 
$\alpha:= \phi + \mu$. Since $\phi + \mu \leq \psi$, we have $\alpha$ is bounded from the above on $\Omega$. In case $\alpha (a) = \alpha(b)$, consider a smooth functiom $\rho$ on $\omega$, supported and positive on a small neighborhood of $a$, and with $\rho (b) = 0$. For a small $t > 0$,, the function 
$\beta := \phi + (\mu +t\rho)$ is PSH, and upper-bounded on $\Omega$, and satisfies $\beta (a) > \beta (b)$.
\qed

\claim 4.4. Relation to ${\bf c}'(\Omega)$ and Bergman spaces|\rm The 
${\bf c}'(\Omega)$ variant of the core was defined in [GHH17]. They denote by 
$PSH'(\Omega)$ the class of all $\alpha\in\PSH(\Omega)$ with Lelong number vanishing identically in $\Omega$. Then $\Omega \setminus {\bf c}'(\Omega)$ is the union of all open subsets $U$ of $\Omega$ such that $\alpha |U$ is {\it strongly} PSH for some upper-bounded $\alpha \in PSH'(\Omega)$. Gallagher et al
[GHH17] show that, for a pseudoconvex $\Omega$ Bergmann space $A^2(\Omega)$ is infinite dimensional if 
${\bf c}'(\Omega) \neq \Omega$ (Remark 8(a)), and if ${\bf c}'(\Omega)$ empty, then $\Omega$ has Bergmann metric (Remark 8(b)). In Remark 8(d) they give an example of such domain (in $\bC^2$)  for which the core ${\bf c}(\Omega)$ is nonempty, specifically a complex line. We provide here a somewhat sronger, analogous example, with 
${\bf c}(\Omega)$ having non-empty interior.
\endclaim

\claim 4.5. Remark|\rm If $\log\mu \geq 0$ is strongly PSH on $\bC^2$, then any domain $\Omega$ as in Corollary 1.6 has ${\bf c}'(\Omega)$ empty, and so has the Bergmann metric. This holds for domains obtained in Section 3, which have 
${\bf c}(\Omega)$ with nonempty interior, with  
$\mu (z, w) := 1 +|\Phi (z,w)|^2$, where $\Phi$ is a Fatou-Bieberbach map. (See Remark 3.4.)
\endclaim

\proof We just apply the argument of [GHH17], Remark 8.4 (d) to more general functions. Since $\psi |\Omega < C$, 
and $\phi + \mu \leq \psi$,
then, letting $\alpha:=\phi +\mu/2 - C$, we obtain $\alpha + \mu /2 < 0 $, in $\Omega$. Since $\alpha < 0$ and 
$\alpha\in PSH(\Omega)$, by [GHH17], Remark 8 (c), $- \log (1 - \alpha)\in PSH'(\Omega)$. Observe that
$$
- \log (1 - \alpha) + \log(1+\mu/2) < 0 \hbox{ , on } \Omega. \leqno (4.5.1)
$$
(This follows from the elementary inequality: 
$$ 
\hbox{if } a+b < 0 \hbox { , and }  b\geq 0 \hbox{ , then }   \log(1-a)+ \log(1+b) < 0 ,
$$ 
applied to numbers $a:= \alpha(z)$, and $b:=\mu(z)/2$, where $z\in\Omega$.)
Note that the negative function (4.5.1) must be in $PSH'(\Omega)$, as the sum of a 
$PSH'$-function and of nonnegative PSH function $\log(1+\mu/2)$. The latter function is {\it strongly} PSH in $\Omega$, as the composition of the increasing convex function $t\to\log(1+e^t/2)$ with the strongly PSH $\log\mu$. Summarizing, (2.5.1) is an upperbounded, $PSH'$-function, and strongly PSH on $\Omega$, and so ${\bf c}'(\Omega)$ is empty.
\qed

\plainsection{A. Appendix}

\claim A.1. Proposition|Let $P(z, w)$ be a nonconstant complex polynomial. Then either

\emph{(i)}  there are at most finitely many complex lines $L_1, L_2, ..., L_s$ in $\bC^2$, that are disjoint from the zero set of $P$; or

\emph{(ii)} there is a complex line $L_0$, such that $P$ is constant on every line $L'$ parallel to $L_0$, and so, every complex line $L$ that is not parallel to $L_0$ must intersect the zero set of $P$.

\endclaim

\proof Let $\cL$ denote the set of all complex lines $L$ in $\bC^2$ such that $P$ is not vanishing on $L$. Then $P|L$ is constant.

{\it Case(a)}. Suppose not all lines in $\cL$ are parallel. We claim that in this case there are at most $d$ lines in $\cL$, where $d$ is the degree of $P$.

Suppose this is false; then we can choose $(d+1)$ distinct lines $L_0, L_1, ..., L_d$ in $\cL$, so that $L_0$ and $L_1$ intersect. Then every $L_i$ intersect $L_0$ or $L_1$, and so $S:=\bigcup_{0\leq i \leq d} L_i$ is connected, and  $P|S = c$ where $c$ is a constant. 
Consider any point $p\in \bC^2 \setminus S$. Then there is a complex line $L^*$ through $p$ that intersects evary line $L_i$, for $i=0, 1 , ..., d$, but does not pass through any of their intersection points. (We skip easy details here.) Then $P|L$ has the same value $c$ at each of the $(d+1)$ distinct intersection points of $L^*$ with $S$. Since degree of the polynomial $P|L^*$ is less than $(d+1)$, $P|L^*$ must be constant, equal $c$. We obtain $P$ is constant on $\bC^2$ contrary to the assumptions.

{\it Case(b)}. Suppose $\cL$ is a set of parallel lines, and there are more then $d$ of them. Then $P$ is constant on every complex line parallel to them, and every complex line transversal to $\cL$ intersects the zero set of $P$.

To see this, consider $(d+1)$ distinct lines, say $L_0, L_1,..., L_d$ in $\cL$. Choose a complex coordinate system in $\cL$, so that $L_i = \{(z,w)\in \bC^2 : w_i=b_i \}$, with $b_0, b_1,...,b_d$ distinct constants. Consider, for any two constants $a, a'$ the polynomial $P(a, w) - P(a', w)$. Since it has $(d+1)$ distinct roots (i.e. $b_j$'s), it is identically zero. Hence, $P$ is constant on every line parallel to $L_0$, and so, $P$ being nonconstant, for every complex line $L$ trannsversal to $L_0$, $P|L$ is a nonconstant polynomial and so must have a zero.

Clearly cases (a) and (b) together imply the proposition.
\qed

\plainsection{References}

\parskip 1.5pt plus 0.5pt minus 0.5pt

\Bibitem[AyS14]&Aytuna, A., Sadullaev, A.:& Parabolic Stein manifolds.& Math. Ann. {\bf 114} (2014), 86-109&
 
\Bibitem[BeT82]&Bedford, E.,Taylor, B. A.:& A new capacity for plurisubharmpnic functions.& Acta mathematica. {\bf 149} (1992), 1-40&
 
\Bibitem[Dem85]&Demailly, J.-P.:& Mesures de Monge-Ampère et caractérisation géométrique des variétés algébriques affines.& Mémoires de la S. M. F. 2e série. {\bf 19} (1985)&
 
\Bibitem[HST17]&Harz, T., Shcherbina, N., Tomassini, G.:& On defining functions and cores for unbounded
domains I.& Math. Z. {\bf 286} (2017), 987-1002&

\Bibitem[HST20]&Harz, T., Shcherbina, N., Tomassini, G.:& On defining functions and cores for unbounded
domains II.& J. Geom. Anal. {\bf 30} (2020), 2293-2325&

\Bibitem[HST21]&Harz, T., Shcherbina, N., Tomassini, G.:& On defining functions and cores for unbounded
domains III.& Mat. Sb. {\bf 212}, 6 (2021), 126-156&

\Bibitem[HST12]&Harz, T., Shcherbina, N., Tomassini, G.:& Wermer type sets and extensions of CR functions.& Indiana Univ. Math. J. {\bf 62} (2012), 431 -459&

\Bibitem[GHH17]&Gallagher, A.-K., Harz, T.,Herbort, G.:& On the dimension of the Bergman space for some unbounded domains .& J. Geom. Anal.{\bf 27} (2017), 1435-1444&

\Bibitem[Glo98]&Globevnik, J.:& On Fatou-Bieberbach domains.& Math. Z. {\bf 229}, 1, (1998), 91-106&

\Bibitem[KeR81]&Kerzman, N., Rosay, J.-P.:&Fonctios plurisousharmoniques d'exhaustions bornees et domaines taut.& Math. Ann. {\bf 257} (1981), 171-184&

\Bibitem[Kol05]&Kolodziej, S.:&The complex Monge-Ampere equation and pluripotential theory.&Mem. of the A.M.S. {\bf 173}
(2005), No. 840&

\Bibitem [MST18]&Mongodi, S., Slodkowski, Z., Tomassini, G.:& Weakly complete complex surfaces.& Indiana  Univ. Math. J. {\bf 67} (2018), 899-935&

\Bibitem[PoS19]&Poletsky, E. A., Shcherbina, N.:& Plurisubharmonically separable complex manifolds. & Proc. Amer. Math. Soc. {\bf 147} (2019), 2413 - 2424&

\Bibitem[Ric68]&Richberg, R.:& Stetige streng pseudokonvexe Funktionen.& Math. Ann. {\bf 175} (1968), 257--286&

\Bibitem[SlT04]&Slodkowski, Z., Tomassini, G.:& Minimal kernels of weakly complete spaces. & J. Funct. Anal. {\bf 210} (2004), 125--147&

\Bibitem[Slo19]&Slodkowski, Z.:&Pseudoconcave decompositions in complex manifolds.& Contemp. Math. {\bf 735} (2019), 239-259&

\Bibitem[Smi86]&Smith, P.:& Smoothing plurisubharmonic functions on complex spaces. & Math.
Ann.\ {\bf 273} (1986), 397--413&

\vskip3mm\noindent
Zbigniew Slodkowski
Department of Mathematics, University of Illinois at Chicago\\
851 South Morgan Street, Chicago, IL 60607, USA\\
\emph{e-mail}\/: zbigniew@uic.edu

\end{document}